\documentclass[5p]{elsarticle}

\usepackage[utf8]{inputenc}
\usepackage[T1]{fontenc}
\usepackage[english]{babel}

\input{preamble}
\makeatletter
\def\ps@pprintTitle{%
	\let\@oddhead\@empty
	\let\@evenhead\@empty
	\def\@oddfoot{}%
	\let\@evenfoot\@oddfoot}
\makeatother

\begin{document}
	\begin{frontmatter}
		\title{The Undirected Two Disjoint Shortest Paths Problem}

%
		\author{Marinus Gottschau\fnref{tum}}
		\ead{marinus.gottschau@tum.de}
		\author{Marcus Kaiser\fnref{tum}}
		\ead{marcus.kaiser@tum.de}
		\author{Clara Waldmann\fnref{tum}}
		\ead{clara.waldmann@tum.de}
		\fntext[tum]{Department of Mathematics, Technische Universität München.}
		\begin{abstract}

The $k$~disjoint shortest paths problem (\dspp[k]) on a graph with $k$~source-sink pairs~$(s_i, t_i)$ asks for the existence of $k$~pairwise edge- or vertex-disjoint shortest \stpaths{i}.
It is known to be NP-complete if $k$ is part of the input.
Restricting to \dspp[2] with strictly positive lengths, it becomes solvable in polynomial time.
We extend this result by allowing zero edge lengths and give a polynomial time algorithm based on dynamic programming for \dspp[2] on undirected graphs with non-negative edge lengths.

		\end{abstract}
		\begin{keyword}
			disjoint paths\sep disjoint shortest paths\sep dynamic programming\sep mixed graphs
		\end{keyword}
	\end{frontmatter}
%

\section{Introduction}

Due to many practical applications, e.g., in communication networks, the $k$~disjoint paths problem (\dpp[k]) is a well studied problem in the literature.
The input of the problem is an undirected graph~$G = (V, E)$ as well as $k$~pairs of vertices $(s_i, t_i)\in V^2$ for~$i \in [k] := \set{1, \ldots, k}$ and the task is to decide whether there exist $k$~paths~$P_1, \ldots, P_k$ such that~$P_i$ is an \stpath{i} and all paths are pairwise disjoint.
Here, disjoint can either mean vertex-disjoint or edge-disjoint.
\\
The $k$ disjoint shortest path problem (\dspp[k]) is a generalization of the $k$~disjoint paths problem.
The input of the problem is an undirected graph~$G = (V, E)$ with edge lengths $\ell: E \to \Rbb$ and $k$~pairs of vertices~$(s_i, t_i) \in V^2$ for~$i\in [k]$.
But here, all paths~$P_i$ for~$i \in [k]$ are additionally required to be shortest \stpaths{i}.
Note, if~$\ell \equiv 0$, this agrees with \dpp[k].
\\
We shall refer to the versions of the problems in directed graphs by \ddpp[k]~and \ddspp[k].

\subsection{Related Work}

Probably most famously, Menger's theorem \cite{Menger} deals with disjoint paths which gave rise to one of the most fundamental results for network flows: the max-flow-min-cut theorem~\cite{maxflowmincut2,maxflowmincut1}.
Using these results, an application of any flow algorithm solves the \ddpp[k] if~$s_i = s_j$ for all~$i, j \in [k]$ or~$t_i = t_j$ for all~$i, j \in [k]$.
Without restrictions on the input instances, all variants of the discussed problems are NP-complete if $k$ is considered part of the input \cite{undirectededgehardness,undirectedvertexhardness}.
\\
Due to this, a lot of research focuses on the setting where $k$ is considered fixed.
\citet{undirected_DPP} came up with an $\Ocal\big(\card{V}^3\big)$~algorithm for \dpp[k].
\\
In contrast to that, \citet{subgraph_homeomorphism} prove that \ddpp[k] is still NP-hard, even if $k = 2$.
They give an algorithm that solves \ddpp[k] for any fixed $k$ on directed acyclic graphs in polynomial time.
\citet{acyclic_mixed_graphs} then extended the work of \citet{subgraph_homeomorphism} to solve the problem on acyclic mixed graphs, which are graphs that contain arcs and edges where directing any set of edges does not close a directed cycle.
\\
Since \ddspp[k] and \ddpp[k] agree for $\ell \equiv 0$, all hardness results carry over. However, if all edge lengths are strictly positive \citet{directed_disjoint_shortest_paths} give a polynomial time algorithm for \ddspp[2].
Also, for \dspp[2] with strictly positive edge lengths a polynomial time algorithm is due to \citet{DSPP}.
However, the complexity of \dspp[k] on undirected graphs with non-negative edge lengths and constant~$k\geq2$ is unknown. We settle the case $k=2$ in this paper.
\\
Other than restricting the paths to be shortest \stpaths{i}, e.g., Suurballe \cite{suurballe} gave a polynomial time algorithm minimizing the total length, if all arc lengths are non-negative and~$s_i = s_j$, $t_i = t_j$ for all~$i, j \in [k]$.
\citet{bjorklund} came up with a polynomial time algebraic Monte Carlo algorithm for solving \dpp[2] with unit lengths where the total length of the paths is minimized.

\begin{table}[!b]
    \centering
    \begin{footnotesize}
    \setlength{\tabcolsep}{0.27em}
    \begin{tabular}{lcccc}
    \toprule
     & \multicolumn{2}{c}{$\ell \equiv 0$} & \multicolumn{2}{c}{$\ell$ non-negative} \\
    $k$ & \dpp[k] & \ddpp[k] & \dspp[k] & \ddspp[k] \\
    \midrule

    arb.& NP-hard \cite{undirectededgehardness,undirectedvertexhardness}  & NP-hard \cite{undirectededgehardness} & NP-hard \cite{DSPP} & NP-hard \cite{DSPP} \\
    fixed & P \cite{undirected_DPP} & NP-hard \cite{subgraph_homeomorphism} & open ($\ell > 0$) & open ($\ell > 0$)\\
    & & & open ($\ell \geq 0$) & NP-hard  ($\ell \geq 0$) \cite{subgraph_homeomorphism}\\
    $2$ & P \cite{undirected_DPP} & NP-hard \cite{subgraph_homeomorphism} & P ($\ell > 0$) \cite{DSPP} & P ($\ell > 0$) \cite{directed_disjoint_shortest_paths} \\
    & & & P ($\ell \geq 0 $) $^*$ & NP-hard ($\ell \geq 0$) \cite{subgraph_homeomorphism}\\
    \bottomrule
    \end{tabular}
    \end{footnotesize}
    \caption{Complexity of the disjoint paths problem and its variants. \\$^*$ A polynomial time algorithm for the \dspp[2] on undirected graphs with non-negative edge lengths is the main result of this paper.}
\end{table}

\subsection{Our Results}
We give a polynomial time algorithm for \dspp[2] on undirected graphs with non-negative edge lengths.
Combining techniques from \cite{subgraph_homeomorphism} and \cite{directed_disjoint_shortest_paths} enables us to deal with edges of length zero.
We consider the following problem.

\begin{problem}[Undirected Two Edge-Disjoint Shortest Paths Problem]
\label{problem:2dspp}
	\emph{\textbf{Input:}}
	An undirected graph~$G = (V, E)$ with non-negative edge lengths~$\ell: E \to \Rbbp$, a tuple of sources $s \in V^2$, and a tuple of sinks $t \in V^2$

	\emph{\textbf{Task:}}
	Decide whether there exist two edge-disjoint  paths~$P_1$ and~$P_2$ in~$G$ such that~$P_1$ is a shortest \stpath{1} and~$P_2$ is a shortest \stpath{2} w.r.t.\ the edge lengths~$\ell$.
\end{problem}

Our paper is organized as follows.
In~\cref{sec:mixedgraphs}, based on the ideas of \cite{subgraph_homeomorphism}, we give a dynamic algorithm that solves the \dpp[k] in polynomial time on weakly acyclic mixed graphs, which are a generalization of directed acyclic graphs.\\
These results are then used in~\cref{sec:undirectedweightedgraphs} together with a similar approach as in \cite{directed_disjoint_shortest_paths} to solve the undirected \dspp[2] with non-negative edge lengths in polynomial time.

The results of this paper have been obtained independently by Kobayashi and Sako.


\section{Disjoint Paths in Weakly Acyclic Mixed Graphs}
\label{sec:mixedgraphs}

In this section, we give an algorithm that solves \dpp[k] in a generalization of directed acyclic graphs.
We first define mixed graphs, introduce some notations, and state the problem.

A graph~$G = (V, A \cupdot E)$ is a \emph{mixed graph} on the vertex set~$V$ with arc set~$A \subseteq V^2$ and edge set~$E \subseteq \genfrac(){0pt}{1}{V}{2}$.
We define~$\textAE(G) := A \cupdot E$.
The set of ingoing (outgoing) arcs of a set of vertices~$W \subseteq V$ is denoted by~$\delta^-_A(W)$~$\big(\delta^+_A(W)\big)$.

For pairwise disjoint vertex sets~$W_1, \ldots, W_h$, we denote by~$G / \set{W_1, \ldots, W_h}$ the graph that results from~$G$ by contracting~$W_1, \ldots, W_h$ into $h$ vertices.

A (directed) \emph{$\path{u}{w}$~$P$} in~$G$ is a sequence of $h$ arcs and edges~$(\textae_1, \ldots, \textae_h) \in \textAE^h$ such that there exists a sequence of vertices~$(u = v_1, \ldots, v_{h + 1} = w) \in V^{h + 1}$ satisfying either~$\textae_i = (v_i, v_{i + 1})$ or~$\textae_i = \set{v_i, v_{i + 1}}$ for all~$i \in [h]$.
Two paths are \emph{arc/edge-disjoint} (\emph{vertex-disjoint}) if they do not have a common arc or edge (vertex).

Note that a directed acyclic graph induces natural orderings of its vertices.
A linear ordering of the vertices is called a \emph{topological ordering} if, for every arc $(v, w)$, the tail $v$ precedes the head $w$ in the ordering.
An ordering is called a \emph{reverse topological ordering} if its reverse ordering is a topological ordering.

On a ground set~$U$, a \emph{binary relation} $R$ is a subset of~$U^2$.
For~$(u, v) \in R$, we write~$u \mathbin{R} v$.
A relation~$R$ is called \emph{reflexive}, if~$u \mathbin{R} u$ holds for all~$u \in U$.
For two binary relations~$R, S \subseteq U^2$, the \emph{composition}~$S \circ R$ is defined by~$\set[(u, w) \in U^2]{\exists v \in U: u \mathbin{R} v \land v \mathbin{S} w }$.
Note that~$\circ$ is an associative operator.

We consider the following problem for fixed~$k$.
\begin{problem}[Mixed $k$~Arc/Edge-Disjoint Paths Problem]
\label{problem:mixed}
    \emph{\textbf{Input:}}
    A mixed graph~$G=(V, \textAE)$, a $k$-tuple of sources $s \in V^k$, and a $k$-tuple of sinks~$t \in V^k$

    \emph{\textbf{Task:}}
    Decide whether there exist $k$~pairwise arc/edge-disjoint paths~$P_1, \ldots, P_k$ in~$G$ such that~$P_i$ is an \stpath{i}, for all~$i \in [k]$.
\end{problem}
We give an algorithm that solves this problem on a class of mixed graphs, that generalize directed acyclic graphs:
\begin{definition}[Weakly Acyclic Mixed Graphs]
    We call a mixed graph~$G = (V, A \cupdot E)$ \emph{weakly acyclic} if the contraction of all edges~$E$ yields a directed acyclic graph without loops.
\end{definition}

Note that a weakly acyclic mixed graph can contain (undirected) cycles in its edge set.
\\
For a mixed graph~$G = (V, \textAE)$, we use the following notation in order to discuss the existence of disjoint paths.

\begin{definition}[Arc/Edge-Disjoint Paths Relation]
    For~$k \in \Nbb$, we define the binary relation~$\edgedisjointpaths{\textAE}$ on the set~$V^k$ as follows.
    For~$v, w \in V^k$, we have~$v \edgedisjointpaths{\textAE} w$ if there exist pairwise arc/edge-disjoint \paths{v_i}{w_i} for all~$i \in [k]$ in~$\textAE$.
    We will also write~$\edgedisjointpaths{G}$ short for~$\edgedisjointpaths{\textAE(G)}$.
\end{definition}

Since paths of length zero are allowed, the relation~$\edgedisjointpaths{\textAE}$ is reflexive.
In general, it is not transitive.
When considering two relations based on two disjoint sets of arcs and edges, however, these two act in a transitive manner.
In that case, the respective underlying arc/edge-disjoint paths from both relations can be concatenated.
The resulting arc/edge-disjoint paths correspond to an element in the composition of the two relations.

\begin{observation}[Partial Transitivity]
\label{obs:partialtransitivity}
    For disjoint arc/edge sets~$\textAE_1, \textAE_2 \subseteq \textAE$ and vectors of vertices~$u, v, w \in V^k$, it holds
    \begin{equation*}
        u \edgedisjointpaths{\textAE_1} v \land v \edgedisjointpaths{\textAE_2} w \implies u \edgedisjointpaths{\textAE_1 \cupdot \textAE_2} w .
    \end{equation*}
\end{observation}

\begin{algorithm}
    \setstretch{1.25}
    \medskip
    \KwInput{weakly acyclic mixed graph~$G = (V, A \cupdot E)$}
    \KwOutput{$\edgedisjointpaths{G}$ on~$V^k$}
    \medskip
    Find connected components~$V_1, \ldots, V_h$ of the subgraph~$(V, E)$ sorted according to a topological ordering of~$G / \set{V_1, \ldots, V_h}$\;
    \lFor{$j = 1, \ldots, h$}{
        \nonl \\ \hskip2\skiptext Compute~$\edgedisjointpaths{G[V_j]}$ using an algorithm for \dpp[k]%
    }
    Initialize~$\edgedisjointpaths{}$ to the relation~$\set[(v, v)]{v \in V^k}$\;
    \lFor{$j = 1, \ldots, h$} {
        \nonl \\ \hskip2\skiptext Update~$\edgedisjointpaths{}$ to~$\edgedisjointpaths{G[V_j]} \circ \edgedisjointpaths{\delta^-_A(V_j)} \circ \edgedisjointpaths{}$%
    }
    \Return $\edgedisjointpaths{}$
    \medskip
    \caption{Dynamic Program for \dpp[k] in Weakly Acyclic Mixed Graphs}
    \label{alg:dynamicprogram:edgedisjoint:weaklyacyclic}
\end{algorithm}

This observation is exploited in~\cref{alg:dynamicprogram:edgedisjoint:weaklyacyclic} in order to solve \cref{problem:mixed} for fixed~$k$ for weakly acyclic mixed graphs.
It computes the relation~$\edgedisjointpaths{G}$ in polynomial time by dealing with the edges and arcs in~$G$ separately.
\\
For the undirected components, i.e., the connected components of the subgraph $(V, E)$, it uses an algorithm for edge-disjoint paths in undirected graphs (e.g.,~\cite{undirected_DPP}) to find the relation~$\edgedisjointpaths{}$ on each component.

\pagebreak
Afterwards, dynamic programming is used to compute~$\edgedisjointpaths{}$ on successively larger parts of the mixed graph.
As~$G$ is weakly acyclic, contracting all undirected components results in an acyclic graph.
The algorithm iterates over the components in a topological ordering.
Based on~\cref{obs:partialtransitivity}, previously found arc/edge-disjoint paths are extended alternately by arcs between components and edge-disjoint paths within one component.
This approach is a generalization of the methods presented in~\cite{subgraph_homeomorphism}.
\begin{figure}[t]
    \centering
    \begin{tikzpicture}[x=0.75cm, y=1.1cm, label distance=-2pt]
        \coordinate (V_1)   at (0, 0);
        \coordinate (V_j-1) at (5, 0);
        \coordinate (V_j)   at (9, 0);
        \draw (V_1)   ellipse[x radius=0.8, y radius=1];
        \draw (V_j-1) ellipse[x radius=1.3, y radius=1];
        \draw (V_j)   ellipse[x radius=1.3, y radius=1];
        \node[below right=0.55 and 0.55 of V_1]   {$V_1$};
        \node[below right=0.55 and 0.9 of V_j-1] {$V_{j-1}$};
        \node[below right=0.55 and 0.9 of V_j]   {$V_j$};
        \coordinate[vertex, below=0.05 of V_1, label={above:$v_2$}] (v_2);
        \coordinate[vertex, above=0.4 of v_2, label={above:$v_1$}] (v_1);
        \coordinate[vertex, below left=0.05 and 0.4 of V_j-1] (x_3);
        \coordinate[vertex, above=0.4 of x_3] (x_1);
        \coordinate[vertex, below=0.4 of x_3] (x_2);
        \coordinate[vertex, label={above:$v_3$}] (v_3) at (x_3);
        \coordinate[vertex, below right=0.05 and 0.4 of V_j-1, label={above:$p_3$}] (p_3);
        \coordinate[vertex, above=0.4 of p_3, label={above:$p_1$}] (p_1);
        \coordinate[vertex, label={above:$p_2$}] (p_2) at (x_2);
        \coordinate[vertex, below left=0.05 and 0.4 of V_j, label={above:$q_1$}] (q_1);
        \coordinate[vertex, above=0.4 of q_1, label={above:$q_2$}] (q_2);
        \coordinate[vertex, below=0.4 of q_1, label={[xshift=12pt]above:$q_3 = w_3$}] (q_3);
        \coordinate[vertex, below right=0.05 and 0.4 of V_j, label={above:$w_1$}] (w_1);
        \coordinate[vertex, above=0.4 of w_1, label={above:$w_2$}] (w_2);
        \coordinate[vertex] (w_3) at (q_3);
        \draw[arc] (v_1) -- (x_1) node [pos=0.72, fill=white, sloped, inner sep=2pt] {$\cdots$};
        \draw[arc] (v_2) -- (x_2) node [pos=0.72, fill=white, sloped, inner sep=2pt] {$\cdots$};
        \foreach \i in {1, 2, 3}{
            \draw (x_\i) -- (p_\i);
            \draw[arc] (p_\i) -- (q_\i);
            \draw (q_\i) -- (w_\i);
        }
    \def \r {1.4}
    \def \ys {2.3}
        \draw [decorate, decoration={brace, amplitude=5pt, raise=\r cm}, shorten >=-2pt, shorten <=-2pt] (v_2) -- (p_3) node[black, midway, below, yshift=\ys cm, xshift= 8pt] {$v \edgedisjointpaths[j-1]{\phantom{\delta^-_A}} p$};
        \draw [decorate, decoration={brace, amplitude=5pt, raise=\r cm}, shorten >=-2pt, shorten <=-2pt] (p_3) -- (q_1) node[black, midway, below, yshift=\ys cm, xshift=34pt] {$p \edgedisjointpaths[\phantom{j-1}]{\delta^-_A(V_j)} q \edgedisjointpaths[\phantom{j-1}]{G[V_j]} w$};
        \draw [decorate, decoration={brace, amplitude=5pt, raise=\r cm}, shorten >=-2pt, shorten <=-2pt] (q_1) -- (w_1);
    \end{tikzpicture}
    \caption{In iteration~$j$ of~\cref{alg:dynamicprogram:edgedisjoint:weaklyacyclic}, relation~$\edgedisjointpaths{}^j$ is built by concatenating previously computed paths \big($\edgedisjointpaths{}^{j-1}$\big), pairwise different arcs to the next component \big($\edgedisjointpaths{\delta^-_A(V_j)}$\big), and undirected edge-disjoint paths in the next component \big($\edgedisjointpaths{G[V_j]}$\big).}
    \label{fig:dynamicprogram:edgedisjoint:weaklyacyclic}
\end{figure}
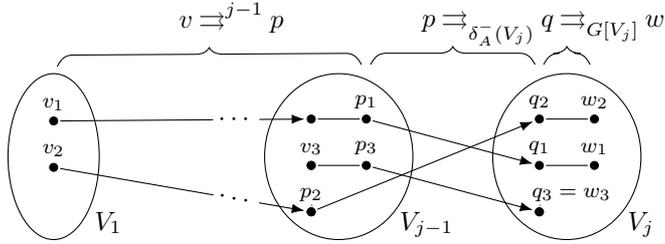

\begin{theorem}[Algorithm\ 1:\ Correctness and Running Time]
    Let~$k \in \Nbb$ be fixed.
    Given a weakly acyclic mixed graph ${G = (V, A \cupdot E)}$, \cref{alg:dynamicprogram:edgedisjoint:weaklyacyclic} computes the relation~$\edgedisjointpaths{G}$ on~$V^k$ in polynomial time.
    \begin{proof}
    Let~$V = \bigcupdot_{j = 1}^h V_j$ be the partition of~$V$ into the vertex sets of the $h$~connected components of~$(V, E)$ as computed by the algorithm.

    For all~$j \in \set{0, \ldots, h}$, let~$\textAE_j$ be the arc and edge set of~$G\big[\bigcup_{l = 1}^j V_l \big]$.
    In particular, $\textAE_0 = \emptyset$ holds true.
    For each~$j \in \set{0, \ldots, h}$, let~$\edgedisjointpaths{}^j$ be the relation~$\edgedisjointpaths{}$ as computed by~\cref{alg:dynamicprogram:edgedisjoint:weaklyacyclic} after the~$j$-th iteration of Line~4.
    In particular, $\edgedisjointpaths{}^0$ is the relation after Line~3.
    In the following, we proof by induction on~$j$ that~$\edgedisjointpaths{}^j$ is equal to~$\edgedisjointpaths{\textAE_j}$.

    After the initialization, this is true for~$j = 0$, as~$\textAE_0$ contains no arcs or edges.
    Consider an iteration~$j \in [h]$ and assume that the claim was true after the previous iteration.

    \enquote{$\subseteq$}:
    Let~$v, w \in V^k$ such that~$v \edgedisjointpaths{}^j w$.
    There exist~$p, q \in V^k$ such that~$v \edgedisjointpaths{}^{j-1} p \edgedisjointpaths{\delta^-_A(V_j)} q \edgedisjointpaths{G[V_j]} w$.
    Using the induction hypothesis, we know~$v \edgedisjointpaths{\textAE_{j-1}} p$.
    Since the arc and edge sets in the three relations are pairwise disjoint, \cref{obs:partialtransitivity} yields~$v \edgedisjointpaths{\textAE_{j}} w$.

    \enquote{$\supseteq$}:
    Let~$v, w \in V^k$ with~$v \edgedisjointpaths{\textAE_j} w$, and~$P_i, i \in [k]$ be arc/edge-disjoint \paths{v_i}{w_i} in~$\textAE_j$.
    Let~$q_i \in V$ be the first vertex on~$P_i$ in~$V_i$ and~$p_i$ be its predecessor if they exist, otherwise set them to~$w_i$ and~$q_i$, respectively.
    As~$G$ is weakly acyclic, we have~$v \edgedisjointpaths{\textAE_{j-1}} p$, $p \edgedisjointpaths{\delta^-_A(V_j)} q$, as well as~$q \edgedisjointpaths{G[V_j]} w$.
    It follows from the induction hypothesis that~$v \edgedisjointpaths{}^j w$.

    The connected components of~$(G, E)$ and their topological ordering in~$G / \set{V_1, \ldots, V_h}$ can be computed in polynomial time.
    Finding edge-disjoint paths in the undirected components can also be done efficiently (e.g.,~\cite{undirected_DPP}).
    A binary relation on~$V^k$ contains at most~$\card{V}^{2 k}$ elements and composing two of them can be done in time polynomial in their sizes.
    Hence, \cref{alg:dynamicprogram:edgedisjoint:weaklyacyclic} runs in time polynomial in the size of the input if~$k$ is fixed.
    \end{proof}
\end{theorem}

In many settings, the problem of finding arc/edge-disjoint paths can be reduced to finding vertex-disjoint paths.
Observe that arc/edge-disjoint paths in a graph correspond to the vertex-disjoint paths in its line graph and an appropriate notion of a line graph can be defined for mixed graphs as well.
\\
For directed graphs, there is a generic reduction from vertex-disjoint to arc-disjoint instances based on splitting vertices.
This reduction, however, cannot be applied to undirected or mixed graphs.
Yet, \cref{alg:dynamicprogram:edgedisjoint:weaklyacyclic} can be modified slightly as follows to compute vertex-disjoint paths.
An algorithm for the undirected vertex-disjoint path problem is used in Line~2.
Only vectors with pairwise different elements are included in the initial relation in Line~3.
Finally in Line~4, tuples~$v, w \in V^k$ are related only if their sets of endpoints~$\set{v_i, w_i}, i \in [k]$ are pairwise disjoint.


\section{Undirected Disjoint Shortest Paths}
\label{sec:undirectedweightedgraphs}

In this section, we study \cref{problem:2dspp} on undirected graphs with non-negative edge lengths.
We first transform the undirected graph~$G$ into a mixed graph and then use the results of the previous section to solve the transformed instance.

\subsection{From Shortest to Directed Paths}

Let an instance of~\cref{problem:2dspp} be given by an undirected graph~$G = (V, E)$, non-negative edge lengths~$\ell: E \to \Rbbp$, and~$s, t \in V^2$.
We are going to transform the graph~$G$ into a mixed graph such that the shortest source-sink-paths in~$G$ correspond to directed source-sink-paths in the resulting mixed graph.

\begin{figure}[h]
	\begin{subfigure}{0.45\columnwidth}
		\centering
		\begin{tikzpicture}
			\node[vertex] (v) [label={left:$v$}] {};
			\node[vertex] (w) [right of=v, label={right:$w$}] {};
			\draw (v) edge (w);
		\end{tikzpicture}
	\end{subfigure}
	\hfill
	$\rightsquigarrow$
	\hfill
	\begin{subfigure}{0.45\columnwidth}
		\centering
		\begin{tikzpicture}
			\node[vertex] (v)  [label={left:$v$}] {};
			\node[vertex] (z-) [above right of=v] {};
			\node[vertex] (z+) [below right of=v] {};
			\node[vertex] (w)  [above right of=z+, label={right:$w$}] {};
			\draw[arc]
				(v)  edge (z-)
				(w)  edge (z-)
				(z-) edge (z+)
				(z+) edge (v)
				(z+) edge (w);
		\end{tikzpicture}
	\end{subfigure}
	\caption{Gadget for resolving conflicts during the orientation of an edge $\set{v,w}$ induced by $d_1$ and $d_2$.}
	\label{fig:opposinggadget}
\end{figure}
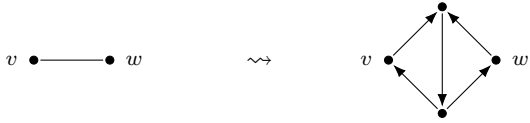

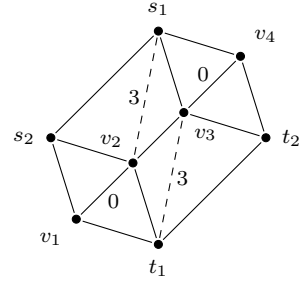
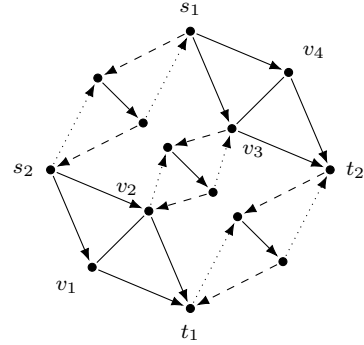
\begin{figure}[t]
	\begin{subfigure}{\linewidth}
		\centering
		\begin{tikzpicture}[font=\footnotesize, node distance=2cm]
			\node[vertex] (s1) [                   label={above:$s_1$}] {};
			\node[vertex] (s2) [below left of=s1,  label={left:$s_2$}] {};
			\node[vertex] (t1) [below right of=s2, label={below:$t_1$}] {};
			\node[vertex] (t2) [below right of=s1, label={right:$t_2$}] {};
			\coordinate (Middles1t2) at ($(s1)!0.5!(t2)$);
			\coordinate (Middles2t1) at ($(s2)!0.5!(t1)$);
			\node[vertex] (2) at ($(Middles2t1) + (45:1.5em)$) [label={100:$v_2$}] {};
			\node[vertex] (5) at ($(Middles2t1) - (45:1.5em)$) [label={below left:$v_1$}] {};
			\node[vertex] (3) at ($(Middles1t2) - (45:1.5em)$) [label={-80:$v_3$}] {};
			\node[vertex] (8) at ($(Middles1t2) + (45:1.5em)$) [label={above right:$v_4$}] {};
			\draw[-]
				(s2) edge (2)
				     edge (5)
				     edge (s1)
				(2) edge (3)
				    edge node[below right, inner sep=1pt]{0} (5)
				    edge (t1)
				(5) edge (t1)
				(3) edge (s1)
				    edge node[above left, inner sep=1pt]{0} (8)
				    edge (t2)
				(s1) edge (8)
				(8) edge (t2)
				(t2) edge (t1)
				;
			\draw[dashed]
				(s1) edge node[left, inner sep=2pt]{3} (2)
				(t1) edge node[right, inner sep=2pt]{3} (3)
				;
		\end{tikzpicture}
		\caption{Example: edges without label have length 1, solid edges are in~$E_1 \cup E_2$}
		\label{fig:ex:usededges}
	\end{subfigure}
	\begin{subfigure}{\linewidth}
		\centering
		\begin{tikzpicture}[font=\footnotesize,node distance=2.6cm]
			\node[vertex] (s1) [                   label={above:$s_1$}] {};
			\node[vertex] (s2) [below left of=s1,  label={left:$s_2$}] {};
			\node[vertex] (t1) [below right of=s2, label={below:$t_1$}] {};
			\node[vertex] (t2) [below right of=s1, label={right:$t_2$}] {};
			\coordinate (Middles1t2) at ($(s1)!0.5!(t2)$);
			\coordinate (Middles2t1) at ($(s2)!0.5!(t1)$);
			\node[vertex] (2) at ($(Middles2t1) + (45:1.5em)$) [label={100:$v_2$}] {};
			\node[vertex] (5) at ($(Middles2t1) - (45:1.5em)$) [label={below left:$v_1$}] {};
			\node[vertex] (3) at ($(Middles1t2) - (45:1.5em)$) [label={-80:$v_3$}] {};
			\node[vertex] (8) at ($(Middles1t2) + (45:1.5em)$) [label={above right:$v_4$}] {};
			\coordinate (Middle23) at ($(2)!0.5!(3)$);
			\node[vertex] (z-23) at ($(Middle23) - (-45:1.2em)$) {} ;
			\node[vertex] (z+23) at ($(Middle23) + (-45:1.2em)$) {};
			\coordinate (Middlet1t2) at ($(t1)!0.5!(t2)$);
			\node[vertex] (z-t1t2) at ($(Middlet1t2) - (-45:1.2em)$) {} ;
			\node[vertex] (z+t1t2) at ($(Middlet1t2) + (-45:1.2em)$) {} ;
			\coordinate (Middles1s2) at ($(s1)!0.5!(s2)$);
			\node[vertex] (z-s1s2) at ($(Middles1s2) - (-45:1.2em)$) {} ;
			\node[vertex] (z+s1s2) at ($(Middles1s2) + (-45:1.2em)$) {} ;
			\draw[arc]
				(s2) edge (2)
					 edge (5)
				(2)	edge (t1)
				(5) edge (t1)
				(s1) edge (3)
				(3)	edge (t2)
				(s1) edge (8)
				(8) edge (t2)
				(z-23) edge (z+23)
				(z-t1t2) edge (z+t1t2)
				(z-s1s2) edge (z+s1s2)
			;
			\draw[arc, dotted]
				(s2)	 edge (z-s1s2)
				(z+s1s2) edge (s1)
				(t1)	 edge (z-t1t2)
				(z+t1t2) edge (t2)
				(2) edge (z-23)
				(z+23) edge (3)
			;
			\draw[arc, dashed]
				(s1)	 edge (z-s1s2)
				(z+s1s2) edge (s2)
				(t2)	 edge (z-t1t2)
				(z+t1t2) edge (t1)
				(3) edge (z-23)
				(z+23) edge (2)
			;
			\draw[-]
				(2) 	edge (5)
				(3)	edge (8)
			;
		\end{tikzpicture}
		\caption{Partially oriented expansion of~(a): solid arcs are in~$A_1 \cap A_2$, dashed arcs are in~$A_1 \setminus A_2$, dotted arcs are in~$A_2 \setminus A_1$}
		\label{fig:ex:gadget}
	\end{subfigure}
	\caption{Exemplary construction of partially oriented expansion}
	\label{fig:ex}
\end{figure}

Since we are interested in shortest $s_1$--$t_1$- and \stpaths{2}, we consider the shortest path networks rooted at~$s_1$ and~$s_2$.
For~$i \in [2]$, we define the distance function~$d_i: V \to \Rbbp$ induced by~$\ell$ w.r.t.\ $s_i$ by~$d_i(v) := \min_{\path{s_i}{v}\,P} \sum_{e \in P} \ell(e)$.
The \emph{shortest path network} rooted at $s_i$ is given by the set
\begin{equation*}
	E_i := \set[\set{v, w} \in E]{\ell(\set{v, w}) = \abs{d_i(v) - d_i(w)}} .
\end{equation*}
See~\cref{fig:ex:usededges} for an example of the sets $E_i$.
\\
The distances~$d_i$ induce an orientation for all edges in~$E_i$ which have a strictly positive length.
We would like to replace an edge~$\set{v, w} \in E$ with~$d_i(v) < d_i(w)$ by the arc~$(v, w)$ (with the same length).
The orientations induced by~$d_1$ and~$d_2$, however, do not have to agree on the set~$E_1 \cap E_2$.
Introducing both arcs would neglect the fact that only one of them can be included in any set of arc/edge-disjoint paths.
We will overcome this by replacing such edges by a standard gadget of directed arcs as depicted in~\cref{fig:opposinggadget}.
\\
Consider the gadget for an edge~$\set{v, w} \in E$.
It contains exactly one \path{v}{w} and one \path{w}{v} corresponding to the two possible orientations of~$\set{v, w}$.
Since both share an arc, only one of two arc/edge-disjoint paths in the transformed graph can use the gadget.
As further both paths consist of three arcs, setting the length of all the arcs in the gadget to~$\tfrac{1}{3} \ell(\set{v, w})$ preserves the distances in the graph.
That way, the distance functions $d_i$ can be extended to the new vertices introduced with gadgets.
\\
For~$i \in [2]$, $A_i$ denotes the set of arcs that result from orienting~$E_i$ w.r.t.\ $d_i$.
More precisely, for~$\set{v, w} \in E_i$ with~$d_i(v) < d_i(w)$ the arc~$(v, w)$ is included into the set~$A_i$ if~$\set{v, w} \in E_1 \triangle E_2$ or the orientation induced by~$d_1$ and~$d_2$ agree.
Otherwise, the arcs of the~\path{v}{w} in the gadget replacing~$\set{v, w}$ are added to~$A_i$.
\\
The induced orientation is only well-defined for edges with strictly positive lengths.
Therefore, the set of edges with length zero~$E_0 := \set[e \in E]{\ell(e) = 0}$ are left undirected and have to be treated in a different manner.

\pagebreak

\begin{definition}[Partially Oriented Expansion]
	Let~$G = (V, E)$ be an undirected graph with non-negative edge lengths~$\ell: E \to \Rbbp$ and~$s \in V^2$.

	The partially oriented expansion of~$G$ w.r.t.\ $\ell$ and~$s$ is the graph~$\smalloverrightharp{G} := (W, E_0 \cup A_1 \cup A_2)$ where
	$W$ is the set of vertices~$V$ augmented with additional vertices introduced with gadgets, and~$E_0$, $A_1$, and~$A_2$ are as defined above.
\end{definition}

The partially oriented expansion of the example from~\cref{fig:ex:usededges} is depicted in~\cref{fig:ex:gadget}.
As we are going to discuss the existence of shortest edge-disjoint paths in~$G$ and the existence of arc/edge-disjoint paths restricted to different arc and edge sets in~$\smalloverrightharp{G}$, the following notation will be useful.

\begin{definition}[Two Disjoint Paths Relations]
\leavevmode\vspace{-1.4em}
	\begin{enumerate}[label=\roman*), wide]
		\item
		Let~$G = (V, E)$ be an undirected graph with non-negative edge lengths~$\ell: E \to \Rbbp$.
		\\
		For~$v, w \in V^2$, we write~$v \edgedisjointshortestpath{\ell}{E} w$ if there exist edge-disjoint shortest \paths{v_i}{w_i} w.r.t.\ $\ell$ for~$i \in [2]$ in~$E$.
		\item
		Let~$G = (V, A \cupdot E)$ be a mixed graph and let~$\textAE_1, \textAE_2$ be two subsets of arcs and edges of $A \cupdot E$.
		\\
		For~$v, w \in V^2$, we write~$v \edgedisjointpathsop[\textAE_1]{\textAE_2} w$ if there exist a \path{v_1}{w_1} in~$\textAE_1$ and a \path{w_2}{v_2} in~$\textAE_2$ which are arc/edge-disjoint.
	\end{enumerate}
\end{definition}

As described above, the distance functions of the original graph~$G$ extend to the vertices of~$\smalloverrightharp{G}$.
For~$i \in [2]$ and~$v \in W$, $d_i(v)$ is the length of a shortest \path{s_i}{v} in~$\smalloverrightharp{G}$.

\begin{lemma}[Paths in the Partially Oriented Expansion]
	\label{lem:partiallyorientedexpansion:equivalence}
	Let~$G = (V, E)$ be an undirected graph with non-negative edge lengths~$\ell: E \to \Rbbp$ and~$s \in V^2$.
	Furthermore, let $\smalloverrightharp{G} = (W, E_0 \cup A_1 \cup A_2)$ be the partially oriented expansion of~$G$ w.r.t.\ $\ell$ and~$s$.

	Then for every~$t \in V^2$, we have~$s \edgedisjointshortestpath{\ell}{E} t$ in~$G$ if and only if~$\genfrac(){0pt}{1}{s_1}{t_2} \edgedisjointpathsop[E_0 \cup A_1]{E_0 \cup A_2} \genfrac(){0pt}{1}{t_1}{s_2}$ in~$\smalloverrightharp{G}$.
	\begin{proof}
		\enquote{$\Rightarrow$}:
		Assume there exist two edge-disjoint shortest \stpaths{i}~$P_i$ in~$E_i$ for~$i \in [2]$.
		Replace each edge with non-zero length in~$P_i$ by the respective oriented arc or path in the respective gadget to obtain~$\smalloverrightharp{P}_i$ in~$E_0 \cup A_i$.
		$\smalloverrightharp{P}_1$~and~$\smalloverrightharp{P}_2$ are arc/edge-disjoint as different edges are replaced by disjoint (sets of) arcs.

		\enquote{$\Leftarrow$}:
		Assume there are arc/edge-disjoint \stpaths{i}~$\smalloverrightharp{P}_i$ in~$E_0 \cup A_i$ for~$i \in [2]$.
		Replace the subpath of~$P_i$ within one gadget with the corresponding edge in~$E_i$.
		The remaining arcs are translated directly to the respective edges in~$E_i$.
		Due to the mentioned equality of distances in~$G$ and~$\smalloverrightharp{G}$ and the fact that~$d_i$ is non-decreasing along arcs in~$\smalloverrightharp{G}$, $P_i$ is a shortest path in~$G$.
		Any path that uses a gadget in~$\smalloverrightharp{G}$, uses its inner arc.
		Therefore, $P_1$ and~$P_2$ inherit being edge-disjoint from~$\smalloverrightharp{P}_1$ and~$\smalloverrightharp{P}_2$.
	\end{proof}
\end{lemma}

\subsection{Disjoint Paths in the Partially Oriented Expansion}

\Cref{lem:partiallyorientedexpansion:equivalence} shows that~$\smalloverrightharp{G}$ captures the shortest paths in~$G$ by using orientation.
We will use the distances, however, to prove the main structural result.
It concern the subgraph of~$\smalloverrightharp{G}$ potentially used by both paths and its weakly connected components, which are its connected components when ignoring the arcs' directions.

\begin{lemma}[Structure of Partially Oriented Expansion]
	\label{lem:partiallyorientedexpansion:structure}
	Let~$G = (V, E)$ be an undirected graph with non-negative edge lengths~$\ell: E \to \Rbbp$ and~$s \in V^2$.
	Furthermore, let $\smalloverrightharp{G} = (W, E_0 \cup A_1 \cup A_2)$ be the partially oriented expansion of~$G$ w.r.t.\ $\ell$ and~$s$.
	Let~$W = \bigcupdot_{j = 1}^h W_j$ be the partition of~$W$ into the vertex sets of the~$h$ weakly connected components of the subgraph~$(W, E_0 \cup (A_1 \cap A_2))$.
\\[1ex]
	Then
	\begin{enumerate}[label=\roman*), wide]
	\item\label{lem:partiallyorientedexpansion:component:weaklyacyclic}
	$\smalloverrightharp{G}[W_j]$ is weakly acyclic for all~$j \in [h]$,
	\item\label{lem:partiallyorientedexpansion:order}
	sorting the components~$W_j, j \in [h]$ in non-decreasing order w.r.t.\ the function $d_1 - d_2$ is a topological ordering of~$(W, A_1) / \set{W_1, \ldots, W_h}$ and a reverse topological ordering of $(W, A_2) / \set{W_1, \ldots, W_h}$, and
	\item\label{lem:partiallyorientedexpansion:component:edges}
	$\smalloverrightharp{G}[W_j]$ contains arcs only from~$A_1 \cap A_2$ and edges only from~$E_0$ for all~$j \in [h]$.
	\end{enumerate}
	\begin{proof}
		\begin{enumerate}[label=\roman*),wide,labelwidth=!,labelindent=0pt]
			\item
			By definition of~$A_1$, we know that~$d_1$ increases strictly along arcs in~$A_1 \cap A_2$.
			Further, $d_1$ is constant on edges in~$E_0$.
			Assume there is~$j \in [h]$ and a (directed) cycle~$C$ in~$\smalloverrightharp{G}[W_j]$ such that there exists~$a \in C \cap A_1 \cap A_2$.
			Along of~$a$ the distance~$d_1$ strictly increases.
			However, $d_1$ cannot decrease along~$C$, which yields a contradiction.
			\item
			Consider the function~ on the vertex set of~$\smalloverrightharp{G}$.
			Based on the common underlying lengths in~$G$ and the definitions of~$A_1$ and~$A_2$, it is strictly increasing along arcs in~$A_1 \setminus A_2$ and strictly decreasing along arcs in~$A_2 \setminus A_1$.
			Opposed to that, it is constant on edges in~$E_0$ as well as along of arcs in~$A_1 \cap A_2$.
			\item
			The function~$d_1 - d_2$ is constant along all arcs~$A_1 \cap A_2$ and edges in~$E_0$.
			Hence, it is constant on each weakly connected component w.r.t.\ those arcs and edges.
			At the same time, the function is not constant along arcs in~$A_1 \triangle A_2$.
			\qedhere
		\end{enumerate}
	\end{proof}
\end{lemma}

This structural result allows to use dynamic programming for solving \cref{problem:mixed} on the partially oriented expansion.
Similar to~\cref{sec:mixedgraphs}, the problem is split into two parts.
First, the two arc/edge-disjoint paths problem on the weakly connected components~$W_1, \ldots, W_h$ of the subgraph~$\left( W, E_0 \cup (A_1 \cap A_2) \right)$ is solved by~\cref{alg:dynamicprogram:edgedisjoint:weaklyacyclic}.
Afterwards, a dynamic program is used to incorporate the results into arc-disjoint paths in~$\smalloverrightharp{G} / \set{W_1, \ldots, W_h}$ to get arc/edge-disjoint paths in~$\smalloverrightharp{G}$.

We know that the two arc/edge-disjoint paths that we are looking for, if they exist, pass through~$\smalloverrightharp{G} / \set{W_1, \ldots, W_h}$ in opposite directions.
In order to accomplish simultaneous construction of both, one of the paths is created backwards.
Apart from that, \cref{alg:dynamicprogram:edgedisjoint:shortest} resembles \cref{alg:dynamicprogram:edgedisjoint:weaklyacyclic}.

\begin{algorithm}
	\setstretch{1.25}
	\medskip
	\KwInput{undirected graph~$G = (V, E)$, non-negative edge lengths~$\ell: E \to \Rbbp$, $s \in V^2$}
	\KwOutput{set of pairs in~$V^2$ that succeed~$s$ w.r.t.\ $\edgedisjointshortestpath{\ell}{E}$}
	\medskip
	Construct~$\smalloverrightharp{G} = (W, E_0 \cup A_1 \cup A_2)$ for~$G$ w.r.t.\ $\ell$ and~$s$\;
	Find weakly connected components~$W_1, \ldots, W_h$ of the subgraph~$\left( W, E_0 \cup (A_1 \cap A_2) \right)$ sorted non-decreasingly w.r.t.\ $d_1 - d_2$\;
	\lFor{$j = 1, \ldots, h$} {
		\nonl \\ \hskip2\skiptext Compute~$\edgedisjointpathsop[{\smalloverrightharp{G}[W_j]}]{\smalloverrightharp{G}[W_j]}$ using~\cref{alg:dynamicprogram:edgedisjoint:weaklyacyclic}%
	}
	Initialize~$\edgedisjointpathsop{}$ to the relation~$\set[(v, v)]{v \in W^2}$\;
	\lFor{$j = 1, \ldots, h$} {
		\nonl \\ \hskip2\skiptext Update~$\edgedisjointpathsop{}$ to~$\edgedisjointpathsop[{\smalloverrightharp{G}[W_j]}]{\smalloverrightharp{G}[W_j]} \circ \edgedisjointpathsop[\delta^-_{A_1}(W_j)]{\delta^+_{A_2}(W_j)} \circ \edgedisjointpathsop{}$%
	}
	\Return $\set[t \in V^2]{\genfrac(){0pt}{1}{s_1}{t_2} \edgedisjointpathsop{} \genfrac(){0pt}{1}{t_1}{s_2}}$
	\medskip
	\caption{Dynamic Program for \dspp[2] with non-negative edge lengths}
	\label{alg:dynamicprogram:edgedisjoint:shortest}
\end{algorithm}

\pagebreak

\begin{theorem}[Algorithm\ 2:\ Correctness and Running Time]
	Given an undirected graph~$G = (V, E)$ with non-negative edge lengths~$\ell: E \to \Rbbp$ and~$s \in V^2$, \cref{alg:dynamicprogram:edgedisjoint:shortest} computes all successors of~$s$ w.r.t.\ $\edgedisjointshortestpath{\ell}{E}$ in polynomial time.
	\begin{proof}
		Let~$W = \bigcupdot_{j = 1}^h W_j$ be the partition of~$W$ into the vertex sets of the~$h$ weakly connected components of the subgraph~$(W, E_0 \cup (A_1 \cap A_2))$ as computed by the algorithm.
		\Cref{lem:partiallyorientedexpansion:structure} \ref{lem:partiallyorientedexpansion:order} shows that the~$W_j$'s are sorted in a topological ordering of~$(W, A_1) / \set{W_1, \ldots, W_h}$ and in a reverse topological ordering of~$(W, A_2) / \set{W_1, \ldots, W_h}$.

		For~$i \in [2]$ and~$j \in \set{0, \ldots, h}$, set~$\textAE^j_i$ to be the arcs of~$A_i$ and edges of~$E_0$ in the induced subgraph~$\smalloverrightharp{G}\big[\bigcup_{l = 1}^j W_l\big]$.
		In particular, we have $\textAE^0_i = \emptyset$.
		For~$j \in [h]$, let~$\edgedisjointpathsop{}^j$ denote the relation~$\edgedisjointpathsop{}$ computed by~\cref{alg:dynamicprogram:edgedisjoint:shortest} after the~$j$-th iteration.
		In particular, $\edgedisjointpathsop{}^0$ is as defined in Line~4.
		We will prove by induction on~$j = 0, \ldots, h$ that~$\edgedisjointpathsop{}^j$ is equal to~$\edgedisjointpathsop[\textAE^j_1]{\textAE^j_2}$.
		The correctness of the algorithm then follows from~\cref{lem:partiallyorientedexpansion:equivalence}.
\\
		The claim holds for~$j = 0$, since~$\textAE^0_1 = \textAE^0_2 = \emptyset$ by definition.
		Consider iteration~$j \in [h]$ and assume that the claim holds for the preceding iteration.

		\enquote{$\subseteq$}:
		Let~$v, w \in W^2$ such that~$v\edgedisjointpathsop{}^j w$.
		Considering Line 5 and using induction hypothesis, there exist~$p, q \in W^2$ with
		\begin{equation*}
			v \edgedisjointpathsop[\textAE^{j-1}_1]{\textAE^{j-1}_2} p \edgedisjointpathsop[\delta^-_{A_1}(W_j)]{\delta^+_{A_2}(W_j)} q \edgedisjointpathsop[{\smalloverrightharp{G}[W_j]}]{\smalloverrightharp{G}[W_j]} w .
		\end{equation*}
		\Cref{lem:partiallyorientedexpansion:structure}~\ref{lem:partiallyorientedexpansion:component:edges} guarantees that the arc and edge sets of the three relations are pairwise disjoint.
		As a result, $v \edgedisjointpathsop[\textAE^j_1]{\textAE^j_2} w$ follows from~\cref{obs:partialtransitivity}.

		\enquote{$\supseteq$}:
		Let~$v, w \in W^2$ such that~$v \edgedisjointpathsop[\textAE^j_1]{\textAE^j_2} w$.
		Thus, there have to be a simple \path{v_1}{w_1}~$P_1$ in~$\textAE^j_1$ and a simple \path{w_2}{v_2}~$P_2$ in~$\textAE^j_2$ that are arc/edge-disjoint.
		Define~$q_1 \in W$ to be the first vertex on~$P_1$ in~$W_j$, if it exists, or~$w_1$.
		Let~$p_1$ be the predecessor of~$q_1$ on~$P_1$ or~$q_1$ if it is the first vertex of~$P_1$.
		Similarly, let~$q_2 \in W$ be the last vertex on~$P_2$ in~$W_j$ or~$w_2$ if it does not exist, and let~$p_2$ be the successor of~$q_2$ or~$q_2$ if~$q_2$ does not have a successor.
		The topological ordering of the~$W_j$'s implies that the subpaths of~$P_1$ and~$P_2$ prove
		\begin{equation*}
			v \edgedisjointpathsop[\textAE^{j-1}_1]{\textAE^{j-1}_2} p \edgedisjointpathsop[\delta^-_{A_1}(W_j)]{\delta^+_{A_2}(W_j)} q \edgedisjointpathsop[{\smalloverrightharp{G}[W_j]}]{\smalloverrightharp{G}[W_j]} w.
		\end{equation*}
		Finally, $v \edgedisjointpathsop{}^j w$ follows by induction hypothesis.

		As for the running time, finding the weakly connected components and sorting them in a topological ordering can be done in polynomial time.
		Computing the relations~$\edgedisjointpathsop[{\smalloverrightharp{G}[W_j]}]{\smalloverrightharp{G}[W_j]}$ also can be done efficiently by virtue of~\cref{alg:dynamicprogram:edgedisjoint:weaklyacyclic}.
		Finally, relations on~$V^2$ have at most~$\card{V}^4$ elements and can be composed efficiently.
		Therefore, the total running time of the algorithm is polynomial in the input size.
	\end{proof}
\end{theorem}

\begin{figure}[t]
	\centering
	\begin{tikzpicture}[x=0.7cm, y=1.2cm, label distance=-2pt]
		\coordinate (W_1)   at (0, 0);
		\coordinate (W_j-1) at (5, 0);
		\coordinate (W_j)   at (9, 0);

		\draw (W_1)   ellipse[x radius=0.8, y radius=0.7];
		\draw (W_j-1) ellipse[x radius=1.2, y radius=0.7];
		\draw (W_j)   ellipse[x radius=1.2, y radius=0.7];

		\node[below right=0.4 and 0.5 of W_1]   {$W_1$};
		\node[below right=0.4 and 0.8 of W_j-1] {$W_{j-1}$};
		\node[below right=0.4 and 0.8 of W_j]   {$W_j$};

		\coordinate[vertex, above=0.15 of W_1, label={above:$v_1$}] (v_1);
		\coordinate[vertex, below=0.25 of W_1, label={above:$v_2$}] (v_2);

		\coordinate[vertex, above left=0.15 and 0.4 of W_j-1] (x_2);

		\coordinate[vertex, above right=0.15 and 0.4 of W_j-1, label={above:$p_2$}] (p_2);
		\coordinate[vertex, below left=0.25 and 0.4 of W_j-1, label={above:$p_1$}] (p_1);

		\coordinate[vertex, above left=0.15 and 0.4 of W_j, label={above:$q_1$}] (q_1);
		\coordinate[vertex, below left=0.25 and 0.4 of W_j, label={[xshift=12pt]above:$q_2 = w_2$}] (w_2);

		\coordinate[vertex, above right=0.15 and 0.4 of W_j, label={above:$w_1$}] (w_1);

		\draw[arc] (v_1) -- (p_1) node [pos=0.72, fill=white, sloped, inner sep=2pt] {$\cdots$};
		\draw[arc] (p_1) -- (q_1);
		\draw (q_1) -- (w_1);

		\draw[arc] (w_2) -- (p_2);
		\draw[arc] (p_2) -- (x_2);
		\draw[arc] (x_2) -- (v_2) node [pos=0.28, fill=white, sloped, inner sep=2pt] {$\cdots$};

		\draw [decorate, decoration={brace, amplitude=5pt, raise=0.8cm}, shorten >=-2pt, shorten <=-2pt] (v_1) -- (p_2) node[black, midway, above, yshift=0.9cm, xshift= 8pt] {$v \edgedisjointpathsop[j-1]{\phantom{\delta^+_{A_2}}} p$};
		\draw [decorate, decoration={brace, amplitude=5pt, raise=0.8cm}, shorten >=-2pt, shorten <=-2pt] (p_2) -- (q_1) node[black, midway, above, yshift=0.9cm, xshift=34pt] {$p \edgedisjointpathsop[\delta^-_{A_1}(W_j)]{\delta^+_{A_2}(W_j)} q \edgedisjointpathsop[{\smalloverrightharp{G}[W_j]}]{\smalloverrightharp{G}[W_j]} w$};
		\draw [decorate, decoration={brace, amplitude=5pt, raise=0.8cm}, shorten >=-2pt, shorten <=-2pt] (q_1) -- (w_1);
	\end{tikzpicture}
	\caption{Iteration~$j$ of~\cref{alg:dynamicprogram:edgedisjoint:shortest}: relation~$\edgedisjointpathsop{}^j$ is built by concatenating already computed paths, pairwise different arcs to the next component, and arc/edge-disjoint paths in the next mixed component}
	\label{fig:dynamicprogram:edgedisjoint:shortest}
\end{figure}
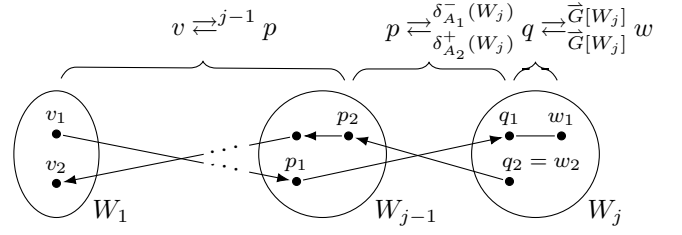

Similar to~\cref{sec:mixedgraphs}, \cref{alg:dynamicprogram:edgedisjoint:shortest} can be adapted to check for the existence of two vertex-disjoint shortest paths.
In that case, the gadget from~\cref{fig:opposinggadget} is not needed anymore, but can be replaced by two opposite arcs.
%
 %
 %
 %
 %

\vspace{-0.5cm}
	\section*{Acknowledgments}
\vspace{-0.35cm}
	This work has been supported by the Alexander von Humboldt Foundation with funds from the German Federal Ministry of Education and Research (BMBF).
	Additionally, we want to thank Jannik Matuschke for his valuable comments and helpful discussions.
\vspace{-0.5cm}
	\bibliographystyle{elsarticle-harv}
	\bibliography{bibliography}
\end{document}